\documentclass[11pt]{article}
\usepackage{geometry}                
\usepackage[parfill]{parskip}    
\usepackage{graphicx}
\usepackage{amsmath, latexsym, amssymb}
\usepackage{epstopdf}
\input xypic
\begin{document}
\title{Manifolds admitting a $\tilde G_2$-structure.}
\author{H\^ong-V\^an L\^e}
\date{}                                           

\maketitle

\newcommand{\R}{{\mathbb R}}
\newcommand{\C}{{\mathbb C}}
\newcommand{\F}{{\mathbb F}}
\newcommand{\Z}{{\mathbb Z}}
\newcommand{\N}{{\mathbb N}}
\newcommand{\Q}{{\mathbb Q}}
\newcommand{\Hq}{{\mathbb H}}

\newcommand{\Aa}{{\mathcal A}}
\newcommand{\Bb}{{\mathcal B}}
\newcommand{\Cc}{{\mathcal C}}    
\newcommand{\Dd}{{\mathcal D}}
\newcommand{\Ee}{{\mathcal E}}
\newcommand{\Ff}{{\mathcal F}}
\newcommand{\Gg}{{\mathcal G}}    
\newcommand{\Hh}{{\mathcal H}}
\newcommand{\Kk}{{\mathcal K}}
\newcommand{\Ii}{{\mathcal I}}
\newcommand{\Jj}{{\mathcal J}}
\newcommand{\Ll}{{\mathcal L}}    
\newcommand{\Mm}{{\mathcal M}}    
\newcommand{\Nn}{{\mathcal N}}
\newcommand{\Oo}{{\mathcal O}}
\newcommand{\Pp}{{\mathcal P}}
\newcommand{\Qq}{{\mathcal Q}}
\newcommand{\Rr}{{\mathcal R}}
\newcommand{\Ss}{{\mathcal S}}
\newcommand{\Tt}{{\mathcal T}}
\newcommand{\Uu}{{\mathcal U}}
\newcommand{\Vv}{{\mathcal V}}
\newcommand{\Ww}{{\mathcal W}}
\newcommand{\Xx}{{\mathcal X}}
\newcommand{\Yy}{{\mathcal Y}}
\newcommand{\Zz}{{\mathcal Z}}

\newcommand{\zt}{{\tilde z}}
\newcommand{\xt}{{\tilde x}}
\newcommand{\Ht}{\widetilde{H}}
\newcommand{\ut}{{\tilde u}}
\newcommand{\Mt}{{\widetilde M}}
\newcommand{\Llt}{{\widetilde{\mathcal L}}}
\newcommand{\yt}{{\tilde y}}
\newcommand{\vt}{{\tilde v}}
\newcommand{\Ppt}{{\widetilde{\mathcal P}}}

\newcommand{\Remark}{{\it Remark}}
\newcommand{\Proof}{{\it Proof}}
\newcommand{\ad}{{\rm ad}}
\newcommand{\Om}{{\Omega}}
\newcommand{\om}{{\omega}}
\newcommand{\eps}{{\varepsilon}}
\newcommand{\Di}{{\rm Diff}}
\newcommand{\Pro}[1]{\noindent {\bf Proposition #1}}
\newcommand{\Thm}[1]{\noindent {\bf Theorem #1}}
\newcommand{\Lem}[1]{\noindent {\bf Lemma #1 }}
\newcommand{\An}[1]{\noindent {\bf Anmerkung #1}}
\newcommand{\Kor}[1]{\noindent {\bf Korollar #1}}
\newcommand{\Satz}[1]{\noindent {\bf Satz #1}}

\newcommand{\gl}{{\frak gl}}
\renewcommand{\o}{{\frak o}}
\newcommand{\so}{{\frak so}}
\renewcommand{\u}{{\frak u}}
\newcommand{\su}{{\frak su}}
\newcommand{\ssl}{{\frak sl}}
\newcommand{\ssp}{{\frak sp}}

\newcommand{\Cinf}{C^{\infty}}
\newcommand{\CS}{{\mathcal{CS}}}
\newcommand{\YM}{{\mathcal{YM}}}
\newcommand{\Jreg}{{\mathcal J}_{\rm reg}}
\newcommand{\Hreg}{{\mathcal H}_{\rm reg}}
\newcommand{\SP}{{\rm SP}}
\newcommand{\im}{{\rm im}}

\newcommand{\inner}[2]{\langle #1, #2\rangle}    
\newcommand{\Inner}[2]{#1\cdot#2}
\def\NABLA#1{{\mathop{\nabla\kern-.5ex\lower1ex\hbox{$#1$}}}}
\def\Nabla#1{\nabla\kern-.5ex{}_#1}

\newcommand{\half}{\scriptstyle\frac{1}{2}}
\newcommand{\p}{{\partial}}
\newcommand{\notsub}{\not\subset}
\newcommand{\iI}{{I}}               
\newcommand{\bI}{{\partial I}}      
\newcommand{\LRA}{\Longrightarrow}
\newcommand{\LLA}{\Longleftarrow}
\newcommand{\lra}{\longrightarrow}
\newcommand{\LLR}{\Longleftrightarrow}
\newcommand{\lla}{\longleftarrow}
\newcommand{\INTO}{\hookrightarrow}

\newcommand{\Sy}{\text{ Diff }_{\om}}
\newcommand{\Ex}{\text{Diff }_{ex}}
\newcommand{\jdef}[1]{{\bf #1}}
\newcommand{\QED}{\hfill$\Box$\medskip}

\newcommand{\UuU}{\Upsilon _{\delta}(H_0) \times \Uu _{\delta} (J_0)}
\newcommand{\bm}{\boldmath}

\begin{abstract}
We find a  necessary and sufficient condition for a compact 7-manifold to admit a $\tilde G_2$-structure.
As a result we find a  sufficient condition for an open 7-manifold to admit a closed 3-form of $\tilde G_2$-type.
 \end{abstract}
\medskip

\section{Introduction}

Recently  a new class of geometries related with stable  forms  has been discovered [Hitchin2000],
[Hitchin2001],  [Witt2005], [Le2006], [LPV2007].  In some cases we can define easily  a necessary
and sufficient condition  for a manifold $M$ to admit  a stable form of  type $\om$ in terms of  topological invariants of $M$, for example if $\om$ is a 3-form  of $G_2$-type [Gray1969].  But in general there is no
method to solve the question how to find a necessary and sufficient condition for a manifold to admit a stable form.  In 
a previous note [Le2006] we have  wrongly  stated   a sufficient condition  for an open manifold to admit a  closed stable 3-form of $\tilde G_2$-type.  We recall that [Bryant1987] a 3-form on
$\R^7$ is called of $\tilde G_2$-type, if it lies on the $Gl(\R^7)$-orbit of a 3-form 
$$\om^3=\theta_1 \wedge \theta_2 \wedge \theta_3 +\alpha_1 \wedge \theta _1 + \alpha_2 \wedge \theta_2 + \alpha_3\wedge \theta _3$$
Here $\alpha_i$ are 2-forms on $V^7$ which can be written as
$$\alpha_1 = y_1\wedge y_2 + y_3\wedge y_4, \: \alpha_2 = y_1 \wedge y_3 - y_2\wedge
y_4, \: \alpha_3 = y_1 \wedge y_4 + y_2 \wedge y_3$$
and $(\theta_1, \theta_2, \theta_3, y_1, y_2, y_3, y_4)$ is an  oriented basis of
$(V^7)^*$.

The group $\tilde G_2$ can be defined as the  isotropy group of $\om^3$ under the action
of $Gl(\R^7)$.  Bryant proved that  [Bryant1987] $\tilde G_2$ coincides with the automorphism
group of  the split octonians.

 In this note we prove the following

\medskip

{\bf Main Theorem}.  {\it Suppose that $M^7$ is a compact 7-manifold. Then  $M^7$ admits
a 3-form of $\tilde G_2$-type, if and only if $M^7$ is orientable  and spinnable.  Equivalently  the
first and second Stiefel-Whitney classes of $M^7$ vanish.  Suppose  that
$M^7$ is an open manifold which admits an embedding to a compact orientable  and spinnable
7-manifold. Then $M^7$ admits a closed 3-form of $\tilde G_2$-type.}

\medskip

\section{Proof of Main Theorem}

Our proof is based on the following simple fact on $\tilde G_2$. 

{\bf 2.1. Lemma.}   {\it We have $\pi_1 (\tilde G_2) = \Z_2$. Hence its maximal compact Lie group
is $SO(4)$.}

\medskip

This Lemma is well-known,  (Bryant mentioned it  but he omitted a proof in [Bryant1987]), but I did not find an explicit  proof of it in  popular lectures on Lie groups, though it
could be given as an exercise. For a hint to a solution of this exercise we refer to [HL1982], p.115,  for an explicit
embedding of $SO(4)$ into $G_2$. The reader can also check that the image of this  group  is also a subgroup
of $\tilde G_2\subset Gl (\R ^7)$.     We shall denote this image by
$SO(4)_{3,4}$.  The Cartan theory on symmetric spaces  implies that $SO(4)_{3,4}$ is a maximal compact Lie
subgroup of $\tilde G_2$.

\medskip

Now let us return to proof of our Main theorem. Clearly if $M^7$ admits a $\tilde G_2$-structure, then it must be orientable and spinnable, since  a maximal compact Lie subgroup $SO(4)_{3,4}$ of
$G_2$ is also a  compact subgroup of $G_2$.
%

\medskip

{\bf 2.2. Lemma.} {\it  Assume that $M^7$ is compact, orientable  and spinnable.  Then $M^7$ admits a
 $\tilde G_2$-structure.}

\medskip

{\it Proof.}  Since $M^7$ is compact, orientable  and spinable, $M^7$ admits a SU(2)-structure [Friedrich1997].
Now it is easy to see that it admits a $SO(4)_{3,4}$-structure, where  $SO(4)_{3,4}$
is a maximal compact Lie subgroup of $G_2$. Hence $M^7$ admits a $\tilde  G_2$-structure.\QED

\medskip


To prove the last statement of the Main Theorem we shall use the  following theorem due to Eliashberg-Mishachev
 to deform the 3-form
$\om^3$ to a closed 3-form $\bar \om^3$ of $\tilde G_2$-type on $M^7$.  

For a subspace $\Rr \subset \Lambda ^p M$  we denote by
$Clo_a \Rr$ a  subspace of the space $Sec \, \Rr$ which consists of closed p-forms
$\om : M \to \Rr$ in the cohomology class $a\in H^p(M)$.

\medskip

{\bf Eliashberg-Mishashev Theorem} [E-M2002,10.2.1]  {\it  Let $M$ be an open manifold,
$a \in H^p (M)$ a fixed cohomology class and  $\Rr$ an open Diff M-invariant subset.
Then the inclusion
$$ Clo _a \Rr \INTO Sec\, \Rr$$
is a homotopy equivalence. In particular,

- any $p$-form $\om : M \to \Rr$ is homotopic  in  $\Rr$ to a closed form $\bar \om$.

- any homotopy of $p$-form $\om_t: M : \to \Rr$ which connects two closed forms
$\om_0, \om_1 \in a$ can be deformed in $\Rr$ into a homotopy of
closed forms $\bar \om_t$ connecting $\om_0$ and $\om_1 \in a$.}
\medskip

Let $\Rr$  be the space
of all 3-forms of $\tilde G_2$-type on  $M = M^7$.  Clearly this space is an open $Diff M^7$-invariant  subset
of $\Lambda ^3 M^7$. Now we apply the  Eliashberg-Mishashev theorem to our  3-form $\om^3$
of $\tilde G_2$-type whose existence  has been proved above.
 Hence $M^7$ admits a closed 3-form $\bar \om^3$ of $\tilde G_2$-type.\QED

\medskip

{\bf 2.3. Remark}. It seems that we can drop the closedness condition in  our Main Theorem and
use the classical  obstruction theory to prove the main Theorem.

\medskip

\section*{Acknowledgement.}  This note is partially supported 
by  grant  of  ASCR Nr IAA100190701.

\medskip

Hong Van Le, Institute of Mathematics, Zitna 25, 11567 Praha 1, hvle@math.cas.cz, \\

\medskip

\begin{thebibliography}{99999}
\bibitem[Adams1996]{Adams1996} {\sc  J.F. Adams}, {\it Lectures on  exceptional Lie groups},
The Chicago University Press, 1996.
\bibitem[Bryant1987] {Bryant1987} {\sc R. Bryant,} {\it Metrics with exceptional
holonomy}, Ann. of Math. (2), 126 (1987), 525-576.
\bibitem[E-M2002]{E-M2002} {\sc Y. Eliashberg and N. Mishachev}, {\it Introduction to the h-Principle}, AMS 2002.
\bibitem[Gray1969]{Gray1969}{\sc A. Gray}, {\it Vector cross products on manifolds},
TAMS 141, (1969), 465-504, (Errata in TAMS 148 (1970), 625).
\bibitem[Friedrich1997]{Fridrich1997} {\sc  Th. Friedrich, I. Kath, A. Moroianu,  U. Semmelmann}, {\it On nearly parallel
$G_2$-manifolds},  Journal Geom. Phys. 23 (1997), 259-286.
\bibitem[HL1982] {HL1982}, {\sc R. Harvey and  H. B. Lawson}, Calibrated geometries, Acta Math. (182), 47-157.
\bibitem[Hitchin2000]{Hitchin2000} {\sc N. Hitchin}, {\it The geometry of
three-forms in 6 and 7 dimensions}, J.D.G. 55 (2000), 547-576.
\bibitem[Hitchin2001]{Hitchin2001} {\sc N. Hitchin}, {\it Stable forms and special metrics}, Contemporean math., (2001), 288,  70-89.
\bibitem[Le2006]{Le2006} {\sc  H. V. Le}, {\it  The existence of symplectic 3-forms on 7-manifolds},
arXiv:math.DG/0603182.
\bibitem[LPV2007]{LPV2007}{\sc H.V.Le, M. Panak and J. Vanzura}, {\it Manifolds admitting stable forms}, in preparation.
\bibitem[Witt2005] {Witt2005} {\sc F. Witt}, {\it  Special metric structures and closed forms}, Ph.D. Thesis , arxiv:math.DG/0502443.
\end{thebibliography}
\end{document}